\documentclass[11pt]{article}
\usepackage[auth-sc]{authblk}
\usepackage{amssymb}
\usepackage{amsthm}
\usepackage{graphicx}
\usepackage{latexsym}
\usepackage{amsmath,amscd}
\usepackage{bm}
\def\part#1{\frac{\partial\phantom{q}}{\partial#1}}

\newenvironment{rmk}{\begin{trivlist}\item[]{\bf Remark:} }
{\end{trivlist}}
\newenvironment{ex}{\begin{trivlist}\item[]{\bf Example:} }
{\end{trivlist}}
\newenvironment{rmks}{\begin{trivlist}\item[]{\bf Remarks:} }
{\end{trivlist}}

\newenvironment{prf}{\begin{trivlist}\item[]{\bf Proof:} }
{\hfill $\Box$ \end{trivlist}}

\newtheorem{thm}{Theorem}

\newtheorem{prp}[thm]{Proposition}

\newcommand{\lie}[1]{\mathfrak{#1}}
\def\End{\mathop{\rm End}\nolimits}

\def\tr{\mathop{\rm tr}\nolimits}

\def\diag{\mathop{\rm diag}\nolimits}
\def\Diff{\mathop{\rm Diff}\nolimits}
\def\SDiff{\mathop{\rm SDiff}\nolimits}

\newcommand{\R}{\mathbf{R}}
\newcommand{\C}{\mathbf{C}}

\newcommand{\PP}{{\mathbf {\rm P}}}

\title{ Higgs bundles and diffeomorphism groups}
 \author{Nigel Hitchin}

 \affil  {Mathematical Institute,
University of Oxford,
Woodstock Road,
Oxford, OX2 6GG}

\begin{document}
 \maketitle
 \thispagestyle{empty}

\let\oldthefootnote\thefootnote
\renewcommand{\thefootnote}{\fnsymbol{footnote}}
\footnotetext{hitchin@maths.ox.ac.uk}
\let\thefootnote\oldthefootnote

 \section{Introduction}
 The prime motivation for this paper is an attempt to find geometric structures on a compact surface $\Sigma$  of genus $g>1$ whose moduli space is described by the higher Teichm\"uller spaces introduced in \cite{Hit1}. These are distinguished components of the space of flat $SL(n,\R)$-connections on $\Sigma$. Each one is diffeomorphic to $\R^N$ where $N=2(n^2-1)(g-1)$ and for $n=2$ it is concretely identifiable as Teichm\"uller space, the moduli space of hyperbolic structures on $\Sigma$. In some ways this goal has already been achieved by others  \cite{GW},\cite{Lab} but in a language far removed from the differential geometry of metrics of constant curvature. We shall here produce a conjectural solution to this question, not for $SL(n,\R)$ for finite $n$, but for what can formally be considered as $SL(\infty,\R)$. The finite rank cases may then be thought of as quantizations of a classical piece of geometry. 
 
 The original discovery of the higher Teichm\"uller spaces used the theory of Higgs bundles, which requires the introduction of a complex structure on $\Sigma$, and an $SO(n)$-bundle with connection $A$. The associated rank $n$ vector bundle $E$ has a holomorphic  orthogonal structure  and we also require a Higgs field $\Phi$, a holomorphic section of the bundle $\End E\otimes K$, where $K$ is the canonical bundle, which is symmetric  with respect to the orthogonal structure. A solution to the equations $F_A+[\Phi,\Phi^*]=0$ defines a flat connection $\nabla_A+\Phi+\Phi^*$ whose holonomy is in $SL(n,\R)$ and a special choice of the holomorphic data $(E,\Phi)$ produces the higher Teichm\"uller space. Our approach is to use this theory replacing gauge groups by  groups of diffeomorphisms. 
 
 Within this context, we may not know what meaning to attach to $SL(\infty,\R)$, but the compact group $SU(n)$ has a well-studied  infinite-dimensional version  $SU(\infty)$ which is defined to be the symplectic diffeomorphisms of the 2-sphere.  We can then introduce $SO(\infty)$ as the subgroup which commutes with reflection in an equator and use Higgs bundle data without worrying about what a flat $SL(\infty,\R)$-connection means. 
 
 Taking the action of $SO(\infty)$ on the sphere, we  replace the principal bundle by an $S^2$-bundle over $\Sigma$  and with this formalism a solution to the Higgs bundle equations analogous to those used to define the higher Teichm\"uller spaces  defines an incomplete  4-dimensional hyperk\"ahler metric  on a disc bundle inside the sphere bundle. It acquires a singularity, a {\it fold}, on the boundary circle bundle, and  extends as a negative-definite hyperk\"ahler metric on the complementary disc bundle, just like the folded K\"ahler metrics in \cite{Bay}. In fact in our case an involution  interchanges the two disc bundles. The three hyperk\"ahler forms  extend as smooth closed 2-forms to the whole sphere bundle.  

Each point of the traditional Teichm\"uller space, in its realization as an $SU(2)$-Higgs bundle,  defines such a hyperk\"ahler metric via the symplectic action of $SU(2)$ on the sphere, and in particular the point corresponding to the given complex structure on $\Sigma$. This is called the canonical Higgs bundle and we call its associated hyperk\"ahler metric   the canonical model.  It has appeared in the mathematics and physics literature several times over the past 30 years, though its global structure has been largely ignored.  The sphere bundle in this case is $\PP(1\oplus K)$ and the fold can be identified with  the unit circle bundle in the cotangent bundle. We conjecture that there exist deformations of this model, where the  boundary of the  disc bundle  embeds nonquadratically in the cotangent bundle. 

There are two existence theorems for the standard Higgs bundle equations: one starts with  the holomorphic data of a stable Higgs bundle $(E,\Phi)$, the other starts with an irreducible flat connection. For the higher Teichm\"uller spaces  the first point of view yields a description of the moduli space as a vector space of holomorphic differentials of different degrees on $\Sigma$. The second characterises the holonomy of the flat connection as positive hyperbolic representations of $\pi_1(\Sigma)$ in $SL(n,\R)$ \cite{FG}. We expect the existence of infinite-dimensional versions of these for the case of $SU(\infty)$. For the first approach, we may observe that given a folded hyperk\"ahler metric on the disc bundle, then  if $\theta$ is the tautological (complex) 1-form on the cotangent bundle,  integrating $\theta^m$  against the symplectic form along the fibres yields a holomorphic section of $K^m$ over $\Sigma$. These sections should determine the metric uniquely. For the second approach, we expect the boundary data on the fold, embedded in the cotangent bundle, to play the role of the holonomy of the flat connection.  
  
We produce three pieces of evidence for this conjecture. 
 For a circle bundle in the cotangent bundle, the real and imaginary parts of the holomorphic symplectic form $d\theta$ restrict to  the data of a  generic pair of closed 2-forms  on a 3-manifold. The first result, an idea of O.Biquard, says that  if these forms are real analytic, then they determine a {\it local} folded hyperk\"ahler extension.  The proof uses twistor theory techniques. 
 The challenge then is to find global boundary conditions which will  produce an extension to the whole disc bundle.
 
  The second piece of evidence is to consider the examples given by the Higgs bundle version of Teichm\"uller space itself. We show that the fold in this case is the unit cotangent bundle of the hyperbolic metric defined by the quadratic differential $q$, and not just for the canonical model $q=0$. We then  describe the infinitesimal deformation of the corresponding hyperk\"ahler metric obtained by varying $q$ and generalise  the resulting formula to an arbitrary differential of degree $m$. The problem, from the first point of view, is then to extend the first order deformations  for $m>2$ to  genuine ones.

  The final evidence is  a global one, and based on the observation that for the higher Teichm\"uller spaces expressed as a sum of vector spaces of differentials the origin is the unique fixed point under the action $\Phi\mapsto e^{i\phi}\Phi$ and this is a standard embedding of the canonical rank 2 Higgs bundle. 
We prove here analogously that the canonical model  is unique among   $S^1$-invariant hyperk\"ahler metrics satisfying the folded boundary condition.The proof  uses a geometrical reformulation of what has long been known as the  $SU(\infty)$-Toda equation, describing  4-dimensional hyperk\"ahler manifolds with this symmetry.

 The author wishes to thank Olivier Biquard and Robert Bryant for useful exchanges of ideas;   the Academia Sinica in Taipei and the Mathematical Sciences  Center, Tsinghua University, Beijing  for their hospitality while part of this paper was being written, and  ICMAT Madrid and QGM Aarhus for support. 
 
 \section{$SU(\infty)$ and $SO(\infty)$}
 We first explain how the group $\SDiff(S^2)$ of symplectic diffeomorphisms of the standard $2$-form $\omega$ on the 2-sphere may be described as $SU(\infty)$. 
 Its Lie algebra consists of $C^{\infty}(S^2)$ modulo the constant functions  -- the Hamiltonian functions for the symplectic vector fields. Equivalently we may consider functions whose integral against $\omega$  is zero. The group $SU(2)$ acts via its quotient $SO(3)$ and breaks up (the $L^2$-completion of) the algebra into an orthogonal sum of irreducible representations ${\bf 3}+{\bf 5}+{\bf 7}+\dots$. Each irreducible representation of $SO(3)$ occurs  with multiplicity one. 
 
 Now consider the  irreducible representation ${\bf n}$  of $SU(2)$. This is a homomorphism $SU(2)\rightarrow SU(n)$ and the Lie algebra of $SU(n)$ under the restricted adjoint action breaks up as ${\bf 3}+{\bf 5}+\dots + {\bf (2n-1)}$. The analogy with $\SDiff(S^2)$ is clear and justifies the notation $SU(\infty)$, although $\lie{su}(n)$ does not embed as a Lie subalgebra in  $\lie{su}(\infty)$: the Poisson bracket does not restrict to the Lie bracket except on the 3-dimensional component $\lie{su}(2).$ 
 
\begin{rmk} The image of this homomorphism is the distinguished {\it principal} 3-dimensional subgroup for $SU(n)$. Higher Teichm\"uller spaces exist for all split real forms $G^r\subset G^c$ \cite{Hit1} and the principal 3-dimensional subgroup for any simple Lie group plays a fundamental role in this construction. Moreover it is the homomorphism of the split groups $SL(2,\R)\rightarrow G^r$ and the corresponding associated  flat connections which maps ordinary Teichm\"uller space into its higher version. 
\end{rmk}
 
 The group $SU(\infty)$ has several formal attributes in common with a compact Lie group, in particular there are invariant polynomials on the Lie algebra which generalize $p_m(A)=\tr A^m$, namely 
 $$p_m(f)=\int_{S^2}f^m\omega$$
 and for $m=2$ this defines a bi-invariant positive definite metric. 
 
 The real form $SL(n,\R)$ is the fixed point set of an involution on $SL(n,\C)$ and restricting to the maximal compact subgroup $SU(n)$, it  fixes the subgroup $O(n)$, with identity component $SO(n)$. Writing the sphere as $x_1^2+x_2^2+x_3^2=1$, and the symplectic form $\omega=2dx_1\wedge dx_2/x_3$ (the induced area form) we define $SO(\infty)$ to be the  subgroup which commutes with the involution $\tau(x_1,x_2,x_3)=(x_1,x_2,-x_3)$ and acts on the equatorial circle $x_3=0$ by an orientation-preserving diffeomorphism. The Lie algebra of $SO(\infty)$ is then formally the functions $f$ which are odd with respect to the involution: since $\tau^*\omega=-\omega$ the Hamiltonian vector fields defined by $i_X\omega=df$ are then even and commute with $\tau$. Note that the induced action on the equator gives a homomorphism $SO(\infty)\rightarrow \Diff(S^1)$.

\section{Higgs bundles}
We recall here the main features of Higgs bundles for a finite-dimensional Lie group and in particular  the choice which yields the higher Teichm\"uller spaces \cite{Hit0},\cite{Hit1},\cite{Sim},\cite{Cor}. Given a compact Riemann surface $\Sigma$ and a compact Lie group $G$ with complexification $G^c$, we take a principal $G$-bundle with connection $A$. To this we can associate a principal $G^c$-bundle and the $(0,1)$ component of the connection defines on it a holomorphic structure. A Higgs field is  a holomorphic section $\Phi$ of ${\lie g}\otimes K$ where ${\lie g}$ denotes the associated Lie algebra bundle and $K$ the canonical line bundle of $\Sigma$. The Higgs bundle equations are $F_A+[\Phi,\Phi^*]=0$ where $F_A$ is the curvature of the connection $A$ and $x\mapsto -x^*$ is the involution on ${\lie g}\otimes K$ induced by the reduction to the compact real form $G$ of $G^c$. 

Given a solution of these equations, the $G^c$-connection with covariant derivative $\nabla_A+\Phi+\Phi^*$ is flat. Conversely, given a flat connection with holonomy a reductive representation $\pi_1(\Sigma)\rightarrow G^c$ there is a reduction of the structure group of the flat principal $G^c$-bundle to $G$ which is defined by a  $\pi_1(\Sigma)$-equivariant harmonic map from the universal covering $\tilde \Sigma$ to the symmetric space $G^c/G$. Using this reduction, the flat connection may be written in the above form for a solution of the Higgs bundle equations. 
 
 If we want a flat connection which corresponds to a representation into a real form $G^r$ of $G^c$, we must take the $G$-connection to reduce to a maximal compact subgroup $H\subset G^r$ and the Higgs field to lie in ${\lie m}\otimes K$ where ${\lie g}=\lie{h}\oplus\lie{m}$. Thus, for $G^r=SL(n,\R)$ we need an  $SO(n)$ connection $A$, or equivalently a rank $n$ vector bundle $E$ with a holomorphic orthogonal structure and $\Lambda^nE$ trivial together with a Higgs field which is symmetric with respect to this structure. 
 
 Uniformization of a Riemann surface gives a representation $\pi_1(\Sigma)\rightarrow PSL(2,\R)$ and a choice of holomorphic square root $K^{1/2}$ of the canonical bundle defines a lift to $SL(2,\R)$.  The corresponding Higgs bundle consists of $E=K^{-1/2}\oplus K^{1/2}$ with the canonical pairing defining the orthogonal structure: $((u,v),(u,v))=\langle u,v\rangle$, and the nilpotent Higgs field $\Phi(u,v)=(v,0)$ is then symmetric. The $SO(2)=U(1)$-connection $A$ is a connection on $K^{1/2}$ which defines a Levi-Civita connection and the Higgs bundle equation $F_A+[\Phi,\Phi^*]=0$ says that the Gaussian curvature is $-4$.

  This canonical example can be modified by taking a holomorphic section $q$ of $K^2$ and defining $\Phi(u,v)=(v, qu)$. The Higgs bundle equations are again for a Hermitian metric $h$ but can be interpreted as saying that the metric 
  \begin{equation}
  \hat h=q+\left(h+\frac{q\bar q}{h}\right)+\bar q
  \label{hyp}
  \end{equation}
  has curvature $-4$ \cite{Hit0}. Then the $3g-3$-dimensional space of quadratic differentials defines Teichm\"uller space from the Higgs bundle point of view.

  For the higher Teichm\"uller spaces of representations into $SL(n,\R)$ we take for  $n=2m+1$ the vector bundle 
 $$E=K^{-m}\oplus K^{1-m}\oplus \dots \oplus K^{m}$$
 and for $n=2m$ 
 $$E=K^{-(2m-1)/2}\oplus K^{1-(2m-1)/2}\oplus \dots \oplus K^{(2m-1)/2}.$$
The pairing of $K^{\pm \ell}$ or $K^{\pm \ell/2}$ defines an orthogonal structure on $E$ and $\Lambda^nE$ is trivial so it has structure group $SO(n,\C)$.

 The Higgs field must be symmetric with respect to this orthogonal structure. We set:
 
  \begin{equation}
 \Phi=\begin{pmatrix}0 & 1 & 0 &\dots & & 0\\
a_2 & 0 & 1 & \dots & & 0\\
a_3 & a_2 & 0 & 1& \dots & 0\\
\vdots & & &\ddots & &\vdots \\
a_{n-1}& & & & \ddots & 1\\
a_n& a_{n-1} & \dots & a_3 & a_2& 0
\end{pmatrix}
 \label{sym}
 \end{equation}
 where $a_i\in H^0(\Sigma, K^i)$. 
 
 Then the higher Teichm\"uller space is the  space  $H^0(\Sigma,K^2)\oplus H^0(\Sigma,K^3)\oplus \cdots \oplus H^0(\Sigma, K^n)$ of differentials of degree $2$ to $n$. The invariants $\tr \Phi^m\in H^0(\Sigma,K^m)$ are universal polynomials in the entries $a_i$ so that given these differentials we have a natural way of defining a Higgs bundle and correspondingly a flat $SL(n,\R)$-connection which fills out a connected component in the moduli space. Setting $a_i=0$ for $i>2$ gives the embedding of Teichm\"uller space. 
 
 Fock and Goncharov \cite{FG} identified this connected component, from the flat connection point of view, with the moduli space of positive hyperbolic representations of $\pi_1(\Sigma)$ in $SL(n,\R)$ (and did so for all split real forms): each homotopy class is mapped to a matrix with $n$ positive eigenvalues. Their methods are quite distinct from the Higgs bundles which lay behind the original discovery.  
 
 \section{Higgs bundles for $SU(\infty)$}
 If we replace the compact group $G$ in a Higgs bundle by $SU(\infty)$ then it is equivalent to consider instead of a principal bundle with structure group $\SDiff(S^2)$, its associated 2-sphere bundle $p:M^4\rightarrow \Sigma$. A connection is defined by a horizontal distribution -- a rank 2 subbundle $H\subset TM$ which is transverse to the tangent bundle  along the fibres $T_F$. Then a horizontal lift of a vector field on $\Sigma$ is a vector field $X$ which integrates to a diffeomorphism taking fibres to fibres. We want this to preserve a symplectic form along the fibres which is a section $\omega_F$ of $\Lambda^2T^*_F$ and so we require ${\mathcal L}_X\omega_F=0$. If $M$ is a symplectic manifold  and the fibres are symplectic submanifolds then the symplectic orthogonal to $T_F$ is an example of a horizontal distribution which defines such an $\SDiff(S^2)$-connection.  Any $SU(2)$-Higgs bundle defines this structure  by virtue of its symplectic action on the sphere. 
 
 For a reduction to $SO(\infty)$ we need an involution $\tau$ on $M$ which acts trivially on the base $\Sigma$ and on each fibre is equivalent by a symplectic diffeomorphism to the standard reflection $(x_1,x_2,x_3)\mapsto (x_1,x_2,-x_3)$. It must preserve the horizontal distribution. The fixed point set is a circle bundle $p:N^3\rightarrow \Sigma$ with a horizontal distribution which we can think of as the connection associated to the homomorphism $SO(\infty)\rightarrow \Diff (S^1)$.   
 
Consider the  Higgs field for $SU(\infty)$ first as a real object: the Lie algebra bundle ${\lie su(\infty)}$ consists of the bundle of smooth functions along the fibres of integral zero. So if $z=x+iy$ is a local coordinate on $\Sigma$  the Higgs field is $\phi_1dx+\phi_2dy$ where $\phi_1,\phi_2$ are local smooth functions on $M$. Then $\Phi=(\phi_1dx+\phi_2dy)^{1,0}$ is to be interpreted as a global section of $p^*K$. If we want an analogue of the Higgs field for a higher Teichm\"uller space then we want $\phi_1,\phi_2$ to be invariant under the involution, or $\tau^*\Phi= \Phi$.

We now need to interpret the Higgs bundle equations in terms of the geometry of this data. A parallel situation was considered  many years ago \cite{AJS} considering Nahm's equations for volume-preserving diffeomorphisms of a 3-manifold.

Consider first the connection -- the horizontal lift of a vector field on $\Sigma$. Locally $\partial/\partial x, \partial/\partial y$ lift to vector fields
$$\frac{\partial}{\partial x}+A_1,\quad \frac{\partial}{\partial y}+A_2$$
where $A_1,A_2$ are Hamiltonian vector fields along the fibres, depending smoothly on $x,y$. The notation is intended to suggest  connection forms in a local trivialization. Replacing the vector fields by Hamiltonian functions $a_1,a_2$ and using the Poisson bracket the equation $F_A+[\Phi,\Phi^*]=0$ reads 
\begin{equation}
(a_2)_x-(a_1)_y+\{a_1,a_2\}+\{\phi_1,\phi_2\}=0.
\label{higgs1}
\end{equation}
The Cauchy-Riemann equation $\bar\partial_A\Phi=0$ is then
\begin{equation}
2(\phi_1+i\phi_2)_{\bar z}+\{a_1+ia_2,\phi_1+i\phi_2\}=0.
\label{higgs2}
\end{equation}

The horizontal distribution splits the tangent bundle of $M$  as $TM\cong T_F\oplus p^*T\Sigma$ and so the relative symplectic form $\omega_F$, a section of $\Lambda^2T^*_F$, defines a genuine 2-form whose restriction to a fibre is $\omega_F$. Using the local formula above for horizontal vector fields this may be written, with $\omega$ the standard symplectic form on $S^2$ and $x_1,x_2$ local coordinates 
$$\omega+\{a_1,a_2\}dx\wedge dy-d_Fa_1\wedge dx-d_Fa_2\wedge dy$$
 where $d_Fa=a_{x_1}dx_1+a_{x_2}dx_2$.  Substituting from the Higgs bundle equation (\ref{higgs1}) gives
 $$\omega-\{\phi_1,\phi_2\}dx\wedge dy-da_1\wedge dx-da_2\wedge dy.$$
Now $\{\phi_1,\phi_2\}dx\wedge dy$ is a well-defined section of $p^*\Lambda^2T^*\Sigma$ and so we may add it to the above form and define on $M$ the closed 2-form locally expressed as 
\begin{equation}
\omega_1= \omega-da_1\wedge dx-da_2\wedge dy.
\label{w1}
\end{equation}
Then
$$\omega_1^2=-2\omega\wedge((a_2)_x- (a_1)_y) dx\wedge dy-2 d_Fa_1\wedge d_Fa_2\wedge dx\wedge dy$$
and from (\ref{higgs1}) this is equal to $2\{\phi_1,\phi_2\}\omega_F\wedge dx\wedge dy$. Hence it is 
 non-degenerate so long as $\{\phi_1,\phi_2\}$ is non-zero. Moreover, in that case, by construction the symplectic orthogonal of $T_F$ is the horizontal distribution. 

Now consider the complex 1-form, the Higgs field $\Phi=(\phi_1+i\phi_2)dz$ and its exterior derivative 
\begin{equation}
\omega^c=d(\phi_1+i\phi_2)\wedge dz.
\label{wc}
\end{equation}
 Then  
$d(\phi_1+i\phi_2)$ and $dz$ span the $(1,0)$-forms for a complex structure so long as $\omega^c\wedge\bar \omega^c\ne 0$. But 
$$\omega^c\wedge\bar \omega^c=d_F(\phi_1+i\phi_2)\wedge d_F(\phi_1-i\phi_2) \wedge d\bar z\wedge dz=4\{\phi_1,\phi_2\}\omega_F\wedge dx\wedge dy$$
and so once more, if $\{\phi_1,\phi_2\}$ is non-zero we obtain a complex structure and a closed holomorphic 2-form $\omega^c$.

Finally consider 
$$\omega_1\wedge \omega^c=\omega\wedge (\phi_1+i\phi_2)_{\bar z}d\bar z\wedge dz-(da_1\wedge dx+da_2\wedge dy)\wedge d(\phi_1+i\phi_2)\wedge dz.$$

From (\ref{higgs2}) the first term is $-\omega_F\wedge \{a_1+ia_2,\phi_1+i\phi_2\}idx\wedge dy$. 
The second term is 
$$(id_Fa_1\wedge d_F(\phi_1+i\phi_2)-d_Fa_2\wedge d_F(\phi_1+i\phi_2))\wedge dx\wedge dy$$
which can be written as 
$\{ia_1-a_2,\phi_1+i\phi_2\}\wedge \omega_F\wedge dx\wedge dy$ and both terms together yield  $\omega_1\wedge \omega^c=0$.

Writing $\omega^c=\omega_2+i\omega_3$  from $(\omega^c)^2=0$ we have $\omega_2^2=\omega_3^2$ and $\omega_2\wedge \omega_3=0$. From $\omega_1\wedge  \omega^c=0$ we have $\omega_1\wedge\omega_2=0=\omega_1\wedge\omega_3$ and with the other calculations we have 
$$\omega_1^2=\omega_2^2=\omega_3^2=2\{\phi_1,\phi_2\}\omega_F\wedge dx\wedge dy \qquad \omega_1\wedge \omega_2=\omega_2\wedge \omega_3=\omega_3\wedge \omega_1=0.$$
It is a standard fact that, given these equations for closed forms,  raising and lowering indices with the symplectic forms and their inverses one can recapture the metric up to a sign and  integrable complex structures $I,J,K$ which define a quaternionic structure on the tangent bundle. This means that we have a {\it hyperk\"ahler metric} on the open set $\{\phi_1,\phi_2\}\ne 0$. 

\begin{ex} Any Higgs bundle for the group $SU(2)$ or $SO(3)$ defines such a hyperk\"ahler metric. The Poisson bracket $\{\phi_1,\phi_2\}$ is then the moment map for a section of the Lie algebra bundle and since moment maps for $SO(3)$ are height functions on the 2-sphere the degeneracy locus of the hyperk\"ahler metric contains  a metric circle bundle, though there may be other points where $[\Phi,\Phi^*]$ itself vanishes. 
\end{ex}

\begin{rmk}
Locally, there is nothing exceptional about these  metrics. Given a 4-dimensional hyperk\"ahler metric, $\omega^c=\omega_2+i\omega_3$ is a holomorphic 2-form with respect to the complex structure $I$ and one may find local holomorphic coordinates such that $\omega^c=dw\wedge dz$. Then regard $z$ as a projection to a local Riemann surface. The Higgs field $\Phi$ is then $wdz$, and the K\"ahler form $\omega_1$ for complex structure $I$ defines the $SU(\infty)$-connection by the symplectic orthogonal to the tangent bundle of the fibres.
\end{rmk}

\section{Folding}
Within symplectic geometry there exists a notion of {\it folded} structure, based on the geometric notion of folding 
a piece of paper along the $y$-axis: the smooth map $f:\R^2\rightarrow \R^2$ defined by $f(x,y)=(x^2,y)$. The standard symplectic form $\omega=dx\wedge dy$ then pulls back to the degenerate form $2xdx\wedge dy$ which is said to be folded along the line $x=0$. In general, a folded symplectic manifold $M^{2m}$ is defined by a closed 2-form $\omega$ such that $\omega^m$ vanishes transversally on a hypersurface $N$ and $\omega$ restricted as a form to $N$ is of maximal rank $2m-2$. There is a Darboux theorem \cite{CGW} which states that locally around the fold $N$ there are coordinates such that 
$$\omega= xdx\wedge dy+\sum_1^{m-1}dx_i\wedge dy_i.$$
 There is also a K\"ahler version of this  and in \cite{Bay} Baykur proves  the remarkable result that any compact smooth 4-manifold has a folded K\"ahler structure, the two components of $M\setminus N$ being Stein manifolds. However, the metric changes signature from positive-definite to negative-definite on crossing the fold. 
 
 For our 4-manifold above, a 2-sphere bundle over the surface $\Sigma$, we have closed 2-forms $\omega_1,\omega_2,\omega_3$ such that $\omega_1^2=\omega_2^2=\omega_3^2$ and this vanishes when the section of $p^*\Lambda^2T^*\Sigma$ defined by $\{\phi_1,\phi_2\}dx\wedge dy$ vanishes. In the case that  the Higgs field is invariant by an involution $\tau$ (the $SO(\infty)$-connection case), the Poisson bracket $\{\phi_1,\phi_2\}$ is anti-invariant and so each term $\omega_i^2$ vanishes on the hypersurface which is the fixed-point set -- a circle bundle over $\Sigma$. We could proceed to define a folded hyperk\"ahler 4-manifold by using the same definition as in the symplectic case, but there is an issue. 
 
 When $\omega_1^2=0$ at a point $x\in M$ we have the algebraic equations in $\Lambda^2T^*_xM$ 
 $$\omega_1^2=\omega_2^2=\omega_3^2=0 \qquad \omega_1\wedge \omega_2=\omega_2\wedge \omega_3=\omega_3\wedge \omega_1=0.$$
If we assume that the $\omega_i$ are  linearly independent at $x$ then they span a 3-dimensional subspace on which the quadratic form given by the exterior product is zero. Now  $\omega^2=0$ is the  condition for $\omega$ to be decomposable -- geometrically it defines a point in the 4-dimensional Klein quadric parametrizing lines in $\PP(T^*_x)$ -- and so we have a plane in the quadric. But there are two types: an $\alpha$-plane consists of the lines through a point (and so $\omega_i=\beta_i\wedge\varphi$) or a $\beta$-plane consists of lines in a plane (or $i_X\omega_i=0$ for some tangent vector $X\in T_xM$). 

If all three closed forms are folded in the usual sense then we would have $\omega_i=xdx\wedge\alpha_i+\beta_i\wedge \gamma_i$ with $\beta_i\wedge\gamma_i$ non-vanishing when restricted as a form to  the fold $x=0$. From the normal form, at $x=0$ we have $i_{\partial/\partial x}(\beta_i\wedge \gamma_i)=0$ and hence we have a $\beta$-plane. 

This is not, however  the case we are interested in -- where the connection is preserved by the involution and the Higgs field is anti-invariant. In this case from (\ref{w1}) we have $\tau^*\omega_1=-\omega_1$ and from (\ref{wc}), $\tau^*\omega_2=\omega_2, \tau^*\omega_3=\omega_3$. In particular $\omega_1$ must vanish as a form when restricted to  the fold and thus is not a folded symplectic structure according to the above definition. 

So take local coordinates such that $x=\{\phi_1,\phi_2\}$ then because $\omega_2,\omega_3$ are even we must have to order $x$ 
$$\omega_2=xdx\wedge\alpha_2+\beta_2\wedge\gamma_2,\quad  \omega_3=xdx\wedge\alpha_3+\beta_3\wedge\gamma_3$$
and 
$$\omega_1=dx\wedge \alpha_1+x\beta_1\wedge \gamma_1$$
where $i_{\partial/\partial x} (\alpha_i,\beta_i,\gamma_i)=0.$
From $\omega_1\wedge \omega_2=0=\omega_1\wedge \omega_3$ we have $\alpha_1\wedge\beta_2\wedge\gamma_2=0=\alpha_1\wedge\beta_3\wedge\gamma_3$ and so we can take $\gamma_2=\gamma_3=\alpha_1=\varphi$. Then at  $x=0$ the $\omega_i$ are all divisible by $\varphi$ and we have explicitly an $\alpha$-plane. 

Also, $d\omega_1=0$ and so $\omega_1=dx\wedge \varphi+xd\varphi$ to order $x$. Furthermore, since $\omega_1^2=2xdx\wedge \varphi\wedge d\varphi$ and vanishes transversally on $x=0$ we have $\varphi\wedge d\varphi\ne 0$ and so $\varphi$  defines a contact structure on the fold $x=0$.

Relabelling $\beta_2,\beta_3$ as $\eta_1,\eta_2$, the restriction of $\omega_2$ and $\omega_3$ to the fold $N$ gives two closed forms $\eta_1\wedge\varphi,\eta_2\wedge\varphi$ where $\varphi$ is a contact form and $\eta_1\wedge\eta_2\wedge\varphi\ne 0$. This latter condition is equivalent to the linear independence of the 2-forms on $M$ at $x=0$.

\begin{rmk}  Hyperk\"ahler metrics with $\alpha$-type folds have been considered in the physics literature. The simplest uses the Gibbons-Hawking Ansatz 
$$g= V\sum_{i=1}^3dx_i^2+V^{-1}(d\phi+\alpha)^2$$
where $V(x_1,x_2,x_3)$ and $\alpha=a_1dx_1+a_2dx_2+a_3dx_3$ satisfy $dV=\ast d\alpha$.  If we take  $V=1/\vert {\mathbf x}+{\mathbf a}\vert + 1/\vert {\mathbf x}-{\mathbf a}\vert $ we get the complete Eguchi-Hanson metric. With a minus sign $V=1/\vert {\mathbf x}+{\mathbf a}\vert - 1/\vert {\mathbf x}-{\mathbf a}\vert $ there is a fold when $V=0$, a 3-manifold given by a circle bundle over the plane through the origin orthogonal to ${\mathbf a}$. There is an involution here covering ${\mathbf x}\mapsto -{\mathbf x}$.

The more general case where $V=\sum_{i=1}^m \pm 1/\vert {\mathbf x}-{\mathbf a_i}\vert$ is folded but the K\"ahler form which vanishes on the fold varies from point to point. This is like the general $SU(2)$-Higgs bundle rather than those which yield flat $SL(2,\R)$-connections.

Metrics like these  can be utilized to construct non-singular Lorentzian signature 5-manifolds which are asymptotically flat and satisfy the equations of 5-dimensional supergravity \cite{GWar}. 
\end{rmk} 
\section{Teichm\"uller space }
\subsection{The canonical model}
The example which we propose to generalize is the canonical $SU(2)$-Higgs bundle described in Section 3 which gives the hyperbolic metric of curvature $-4$ with the given conformal structure on $\Sigma$:

$$E=K^{-1/2}\oplus K^{1/2}\qquad \Phi= \begin{pmatrix}0 & 1 \\
0& 0
\end{pmatrix}.$$

From the inclusion $SO(3)\subset SU(\infty)$ this defines a folded hyperk\"ahler metric on the 2-sphere bundle $M=\PP(K^{-1/2}\oplus K^{1/2})$. Moreover, since the $SU(2)$-connection reduces to $SO(2)$ and the Higgs field lies in ${\lie{m}}\otimes K$ we have the involution $\tau: M\rightarrow M$ defined by $\tau([u,v])=[\bar v h^{-1/2}, \bar u h^{1/2}]$ where $h$, a section of $K\bar K$, is the hyperbolic metric. If we write $M=\PP(1\oplus K)$ and $w=v/u$ then the fixed point set is $h^{-1}w\bar w=1$, the unit circle bundle in the cotangent bundle. We have $[\Phi,\Phi^*]=\diag(h,-h)$ and $\{\phi_1,\phi_2\}$ is the moment map for $h(i,-i)$ which vanishes only on the circle bundle which is the fixed point set of $\tau$. The corresponding positive-definite hyperk\"ahler metric is thus defined on the unit disc bundle in the cotangent bundle, and has an $\alpha$-fold on the unit circle bundle boundary.

As we described in Section 4, the Higgs field $\Phi$ can, in the $SU(\infty)$-interpretation, be thought of as a section of $p^*K$ or equivalently as a map $f:M\rightarrow T^*\Sigma$ to the total space of the cotangent bundle in which case  we can take $\Phi=f^*\theta$ for $\theta=wdz$ the canonical holomorphic 1-form. In our case $f$ is a diffeomorphism restricted to the open disc bundle. 

The canonical Higgs bundle has the property that the map $\Phi\mapsto e^{i\theta}\Phi$, which preserves the equations, takes the solution to a gauge-equivalent one. It follows that the hyperk\"ahler metric is invariant under the action $w\mapsto e^{i\theta}w$, scalar multiplication by a unit complex number on the fibres of $T^*\Sigma$. 
This metric then fits into the more general result of Feix and Kaledin \cite{Feix},\cite{Kal} that a real analytic K\"ahler metric has a unique $S^1$-invariant hyperk\"ahler extension to a neighbourhood of the zero section of the cotangent bundle, where $\omega_2+i\omega_3$ is the canonical holomorphic symplectic form on the cotangent bundle. This is a local result. The canonical model is thus the unique hyperk\"ahler extension of  a metric on $\Sigma$, and  because  the universal covering  is invariant under $PSL(2,\R)$ this must be a  hyperbolic metric. 

The corresponding $SO(3)$-invariant metric for the extension of the round metric on $S^2$ is the Eguchi-Hanson metric. This is complete and extends to the whole cotangent bundle. The hyperbolic analogue, which we are considering here, is less well-known but certainly appears  in the physics literature in \cite{Geg},\cite{Ped}. It was also used by Donaldson \cite{Don} to investigate the hyperk\"ahler extension of the Weil-Petersson metric on Teichm\"uller space. 

We describe this metric first following \cite{Ped}, but with different coordinates. It will reappear in several forms later on. 
We use the upper half-plane model with coordinate $z=x+iy$ and the metric $(dx^2+dy^2)/y^2$ (for convenience of formulae we take constant curvature $-1$ not $-4$) and $w=u_1+iu_2$ the fibre coordinate on the cotangent bundle, with $\theta=wdz$ the canonical 1-form. Then 
\begin{equation}\omega_1= d(u_3(dx-yd\phi))
\label{om1}
\end{equation}
where $\phi=\arg w$ and $y^2(u_1^2+u_2^2+u_3^2)=1$. This latter expression (put $x_i=yu_i$) is the equation for the $S^2$-fibre, and if we put $\varphi=dx-yd\phi$ we see that $\omega_1=du_3\wedge\varphi+u_3d\varphi$ making  $u_3=0$  the fold. We then have 
$$\omega_2+i\omega_3=dw\wedge dz=w(dx+idy)\wedge \left(\frac{1}{y(1-y^2u_3^2)}(dy+y^3u_3du_3)-id\phi\right)$$
and  to first order in $u_3$ this is 
$$-wy^2u_3du_3\wedge dz+w\left(\frac{dy}{y}-id\phi\right)\wedge (dx-yd\phi)$$
and we observe that $\varphi=dx-yd\phi$ is a contact form on the unit circle bundle.
In fact $d(dx/y)=dx\wedge dy/y^2$ so $d\phi-dx/y$ is a connection form for the circle bundle as a principal $S^1$-bundle. 

Restricted as forms to $u_3=0$ we have $\omega_1=0$ and 
 $$
\omega_2=\frac{1}{y^2}(\sin \phi dx+\cos \phi dy)\wedge (dx-yd\phi),\quad 
\omega_3=\frac{1}{y^2}(-\cos \phi dx+\sin \phi dy)\wedge (dx-yd\phi).
$$
The null foliation of $\omega_2$ (which is the real canonical 2-form on $T^*\Sigma$) is tangential to the geodesic flow: in fact explicitly the two equations $\sin \phi \,dx+\cos \phi \,dy=0, dx-yd\phi=0$ give the hyperbolic geodesics $y=c_1\cos\phi, x=c_1\sin\phi+c_2$. We therefore encounter the hyperbolic metric from the data on the fold through its geodesics. 

\begin{rmk} Note that the forms $\omega_2$ and $\omega_3$  are folded in the usual sense, and $\omega_2$ is the K\"ahler form for the complex structure $J$. With a circle action as above the moment map for $\omega_1$ (which from the above calculation is $-u_3y=-x_3$) is a K\"ahler potential for $\omega_2$, and it follows that the complex structure $J$ is Stein and fits in with Baykur's results \cite{Bay}. By contrast the complex structure $I$ admits the compact holomorphic curve given by the zero section of $K$ and is not Stein. 

\end{rmk}

There is a more general picture which we address next.

\subsection{Quadratic differentials}

The description of Teichm\"uller space in the theory of Higgs bundles uses the following pair:
$$E=K^{-1/2}\oplus K^{1/2}\qquad \Phi= \begin{pmatrix}0 & 1 \\
q& 0
\end{pmatrix}$$
where $q$ is a holomorphic section of $K^2$. The canonical model is the case $q=0$. The metric given by solving the Higgs bundle equations is again a Hermitian metric $h$ on $\Sigma$, respecting  the holomorphic orthogonal structure on $E$ and the Higgs field is symmetric, thus from the inclusion $SO(3)\subset SU(\infty)$ this defines a folded hyperk\"ahler metric on the same $S^2$-bundle over $\Sigma$. The actual Higgs bundle equation is 
\begin{equation}
F=2(1-\vert q\vert^2)\omega
\label{qhiggs}
\end{equation}
where $\omega$ is the volume form of the metric and $F$ the curvature of $K$. 
We shall investigate the geometry of the fold in this case. 

Locally, we write $q=a(z)dz^2$ and the metric as $h=kdzd\bar z$. Then  $(k^{-1/4}dz^{-1/2},k^{1/4}dz^{1/2})$ is a unitary basis and relative to this basis 
$$\Phi= \begin{pmatrix}0 & k^{1/2} \\
ak^{-1/2}& 0
\end{pmatrix}dz=k^{1/2}\begin{pmatrix}0 & 1 \\
0& 0
\end{pmatrix}dz+ak^{-1/2}\begin{pmatrix}0 & 0 \\
1& 0
\end{pmatrix}dz.$$
and then, as a 1-form with values in Hamiltonian functions we have 
\begin{equation}
\Phi=(\phi_1+i\phi_2)dz=\frac{1}{2}\left(k^{1/2}(x_1-ix_2)dz+ak^{-1/2}(x_1+ix_2)dz\right)
\label{phiham}
\end{equation}
where $x_1,x_2$ are standard height functions on the unit sphere with $\{x_1,x_2\}=2x_3$. 

If we regard $\Phi$ as a map to $T^*\Sigma$ coordinatized by  $(w,z)\mapsto wdz$ then 
$$x_1-ix_2=2\frac{k^{1/2}w-ak^{-1/2}\bar w}{k-\vert a\vert^2k^{-1}}.$$
The fold $x_3=0$ is then where $x_1^2+x_2^2=1$ which gives the ellipse in the fibre of $T^*\Sigma$
$$\frac{4}{(k-\vert a\vert^2k^{-1})^2} (\bar aw^2+a\bar w^2 +(k+k^{-1}\vert a\vert^2)w\bar w)=-1.  $$
Inverting the matrix of coefficients of $w^2,\bar w^2, w\bar w$ this is the unit circle bundle in the cotangent bundle for the metric
$$\hat h=adz^2+(k+\vert a\vert^2k^{-1})dzd\bar z+\bar ad\bar z^2.$$
But, as in (\ref{hyp}) the Higgs bundle equations imply that this is a metric of constant curvature $-4$. Hence the data on the fold describes the geodesic flow for the hyperbolic metric determined by the quadratic differential $q=adz^2$. 

\section{Local existence}

We have seen in Section 5 that an $\alpha$-folded hyperk\"ahler 4-manifold with involution induces  a contact  structure on the fold $N^3$ defined by a form $\varphi$ and two closed 2-forms $\eta_1\wedge\varphi,\eta_2\wedge\varphi$ with $\eta_1\wedge\eta_2\wedge\varphi\ne 0$. (Note that there is an ambiguity $\varphi\mapsto f\varphi, \eta_1\mapsto f^{-1}\eta_1+g_1\varphi, \eta_2\mapsto f^{-1}\eta_2+g_2\varphi$ in the definition of the 1-forms). In fact, for real analytic data this is sufficient to find a local folded hyperk\"ahler metric on $N\times (-\epsilon,\epsilon)$ for some interval $(-\epsilon,\epsilon)$. I owe the idea below to Olivier Biquard who dealt with a similar issue, with CR boundary data, in \cite{Biq}.

To find a hyperk\"ahler metric we shall use the twistor construction \cite{HKLR}. This means finding a complex 3-manifold $p:Z\rightarrow \PP^1$ fibring over the projective line together with  a holomorphic section $\varpi$ of $\Lambda^2T_F^*(2)$, where $T_F$ is the tangent bundle along the fibres and the factor 2 indicates the tensor product with $p^*{\mathcal O}(2)$, the line bundle of degree $2$ on $\PP^1$. We also need  an antiholomorphic involution $\sigma$, which covers the antipodal map $\zeta\mapsto -1/\bar\zeta$ on $\PP^1$. This defines a  real structure on $Z$ and associated geometrical objects. Then on a family of real sections with normal bundle ${\mathcal O}(1)\oplus {\mathcal O}(1)$ there exists a hyperk\"ahler metric. 

\begin{thm}  Let $X$ be a real analytic 3-manifold with two  analytic closed forms $\eta_1\wedge\varphi, \eta_2\wedge\varphi$ such that the annihilator of $\varphi$ is a contact distribution and  $\eta_1,\eta_2$ satisfy $\eta_1\wedge\eta_2\wedge \varphi\ne 0$. Then this data defines naturally an $\alpha$-folded  hyperk\"ahler metric with involution $(x,t)\mapsto (x, -t)$ on $X\times (-\epsilon,\epsilon)$. 
\end{thm}

\begin{prf} From the real analyticity we complexify locally to a complex manifold $X^c$.

\noindent 1. Consider $\eta=(\eta_1+i\eta_2)-\zeta^2(\eta_1-i\eta_2)$ where $\zeta\in \C$.  
On $\PP^1\times X^c$ we can extend this to a section of $T^*X^c(2)$ by defining $\tilde \eta =-(\eta_1-i\eta_2)+\tilde\zeta^{2}(\eta_1+i\eta_2)$ for $\tilde\zeta=\zeta^{-1}$. Then $\tilde\eta=\zeta^2\eta$ and the pair define a global section. This is moreover real under the conjugation map $(\eta,\zeta)\mapsto (-\bar\eta,-1/\bar\zeta)$. Since $\eta$ is even in $\zeta$ it is also invariant under the holomorphic involution  $\tau$ on $\PP^1\times X^c$  defined by $(\zeta,x)\mapsto(-\zeta,x)$.

Now $\eta\wedge\varphi$ cannot be zero   for taking the exterior product with $\eta_1-i\eta_2$ we have  $$(\eta_1+i\eta_2)\wedge (\eta_1-i\eta_2)\wedge \varphi=-2i\eta_1\wedge\eta_2\wedge\varphi\ne 0.$$   It is moreover real since $\varphi$ is real.  Its annihilator is a holomorphic  1-dimensional subbundle of $TX^c$ over  $\PP^1\times X^c$ and hence describes a foliation by curves.

Locally there is a well-defined  quotient space of the foliation which is a complex 3-manifold $Z$ with a projection to $\PP^1$ induced from the first factor in $\PP^1\times X^c$. It has a real structure and a holomorphic involution induced by $\tau$ whose fixed point set consists of the fibres over $\zeta=0,\infty$. 

Since the form $\eta\wedge\varphi$ is closed on the three-dimensional manifold $X^c$,  this is the quotient by  its  degeneracy foliation and thus each fibre has an induced  symplectic  structure. We therefore have 
  a  holomorphic section $\varpi$ of $\Lambda^2T_F^*(2)$ on $Z$. There are distinguished holomorphic sections of the fibration, the images of $\PP^1\times \{x\}$. These form a 3-dimensional family but what we need is a 4-dimensional family with normal bundle  ${\mathcal O}(1)\oplus {\mathcal O}(1)$.

 \noindent 2. On $\PP^1\times X^c$ the normal bundle of a section  $\PP^1\times \{x\}$ is trivial and since the subbundle defined by the foliation is isomorphic to ${\mathcal O}(-2)$ the normal bundle $N$ in $Z$ fits in an exact sequence 
 $$0\rightarrow {\mathcal O}(-2)\stackrel{\alpha}\rightarrow {\mathcal O}^3\rightarrow N\rightarrow 0.$$
 But $\varphi$ defines a trivial bundle which annihilates the foliation and is thus a trivial subbundle of $N^*$.  The map $\alpha$ is defined by three sections $(s_1,s_2,s_3)$ of ${\mathcal O}(2)$ but the trivial annihilator means that there is a linear relation. With a change of basis the map is $(t_1,t_2,0)$ and $N$ is ${\mathcal O}\oplus L$ where $L$ is the quotient of ${\mathcal O}^2$ by ${\mathcal O}(-2)$ -- in other words ${\mathcal O}(2)$. Thus $N\cong {\mathcal O}(2)\oplus {\mathcal O}$. In fact each of these sections is preserved by the involution $\tau$ so we could take $t_1=1,t_2=\zeta^2$.

This is not the right normal bundle for a hyperk\"ahler metric but in any case we are expecting the metric to be singular for these sections. However,  Kodaira's deformation theory tells us that, given one curve $\PP^1\subset Z$,  as long as $H^1(\PP^1,N)=0$ there is a smooth family of deformations $M^c$  whose tangent space at a curve is isomorphic to the space of holomorphic sections of the normal bundle $N$. For $N={\mathcal O}(2)\oplus {\mathcal O}$, $H^1(\PP^1,N)$ indeed vanishes and $\dim H^0(\PP^1,N)=1+3=4$ so we have a 4-manifold $M^c$.  The real members of this family which have normal bundle ${\mathcal O}(1)\oplus {\mathcal O}(1)$ will define a  hyperk\"ahler manifold. Note that the 3-dimensional family $X^c$ has tangent space $(a_0, a_1\zeta^2-a_2)\in H^0(\PP^1,{\mathcal O}(2)\oplus {\mathcal O})$, the fixed point subspace under $\tau$. The involution induces one on $M^c$ and the fixed point set is $X^c$. 

The condition for a rank 2 holomorphic vector bundle $E$ over $\PP^1$ with $c_1(E)=0$ to be trivial is $H^0(\PP^1,E(-1))=0$. In our case, since the curves are sections of $Z\rightarrow \PP^1$, the normal bundle is the tangent bundle along the fibres $T_F$ and $c_1(T_F)\cong p^*{\mathcal O}(2)$  so the condition for the holomorphic structure to be ${\mathcal O}(1)\oplus {\mathcal O}(1)$ is $H^0(\PP^1,T^*_F)=0$. There is a determinant line bundle over $M^c$ for the $\bar\partial$-operators on $T^*_F$ with a determinant section which vanishes when $H^0(\PP^1,T^*_F)\ne 0$. Unless the determinant is identically zero this is a divisor and the  smooth 3-manifold $X^c$ is already contained in it, so to show that points in $M^c\backslash X^c$ sufficiently close to $X^c$ have normal bundle  ${\mathcal O}(1)\oplus {\mathcal O}(1)$ we have to show that ${\mathcal O}(2)\oplus {\mathcal O}$  is not the generic case. We need the contact condition to do this. 

 \noindent 3. So suppose for  a contradiction that all deformations of the lines parametrized by $X^c$ have normal bundle ${\mathcal O}(2)\oplus {\mathcal O}$. The ${\mathcal O}(2)$ factor is canonically determined as the maximal destabilizing subbundle. Then the sections of ${\mathcal O}(2)$ define a rank 3 distribution on the 4-manifold $M^c$. Take one line $L$ of the family and a one-parameter family $L_t$ passing through a point $z\in L\subset Z$. Then because $z$ is fixed, the  section of the normal bundle $N$ associated to the infinitesimal variation vanishes at $z$ (for all $L_t$) and so lies in the ${\mathcal O}(2)$ component. Hence the whole curve in $M^c$ is tangential to the distribution. Now consider all lines through $z$. This is a 2-dimensional family whose tangent space  at $L$ consists of  the sections of ${\mathcal O}(2)$ which vanish at $z$ -- the quadratic polynomials with factor $(x-z)$. This surface is everywhere tangential to the distribution. It is locally defined by a 1-form $\gamma$, so take two vector fields $U,V$ in the surface so that $i_U\gamma=i_V\gamma=0$, then $i_{[U,V]}\gamma=0$ and $d\gamma$ vanishes as a form on the surface.  
 
 Now the 2-dimensional subspaces of $H^0(\PP^1,{\mathcal O}(2))$ consisting of sections which vanish at points 
$z\in \PP^1$ give elements which span $\Lambda^2(H^0(\PP^1,{\mathcal O}(2)))^*$. Indeed under the $PSL(2, \C)$-invariant isomorphism with $H^0(\PP^1,{\mathcal O}(2))$ these  are polynomials of the form $a(x-z)^2$ and any quadratic is a sum of squares. It follows that  $d\gamma$ vanishes on the whole distribution so that $\gamma\wedge d\gamma=0$ and the distribution is integrable -- a foliation. 

 The tangent space of $X^c$ is, as remarked above,  $(a_0, a_1\zeta^2-a_2)\in H^0(\PP^1,{\mathcal O}\oplus {\mathcal O}(2))$. So the leaves of this foliation intersect $X^c$ tangential to the distribution $a_0=0$ defined by $\varphi$. But $\varphi\wedge d\varphi\ne 0$ so we have a contradiction to integrability. Thus a  twistor line given by a point in the family near $X^c$ must have normal bundle ${\mathcal O}(1)\oplus {\mathcal O}(1)$.

\noindent 4. The three hyperk\"ahler forms may be considered  as the coefficients of a quadratic polynomial $\omega_{\zeta}=(\omega_2+i\omega_3)+2\zeta\omega_1-\zeta^2(\omega_2-i\omega_3)$ and in the twistor approach this is obtained as follows. Since the normal bundle $T_F$ to a twistor line $L_m$  is ${\mathcal O}(1)\oplus {\mathcal O}(1)$ then the natural map 
\begin{equation}
H^0(L_m,T_F(-1))\otimes H^0(L_m,p^*{\mathcal O}(1))\rightarrow H^0(L_m,T_F)= T_mM^c
\label{tensor}
\end{equation}
is an isomorphism. A tangent vector can then be considered as a linear function in $\zeta$ with coefficients sections of $T_F(-1)$.  The section $\varpi$ of $\Lambda^2T^*_F(2)$ is a skew form on the first factor and so evaluating on a pair of tangent vectors gives a quadratic polynomial in $\zeta$. This defines $\omega_{\zeta}$. When the normal bundle jumps, the homomorphism above is no longer an isomorphism, but we can still define $\omega_{\zeta}$ as we shall see next.

\noindent 5. The argument in paragraph 2 of this proof concerning normal bundles only uses the condition $\varphi\wedge d\varphi\ne 0$ and not the full integrability of the leaves of the foliation. It is just the 1-jet of $\varphi$ which is relevant. This means that the section of $T^*_F$ on a line $L_x$ for $x\in X^c$ does not extend to the first order neighbourhood. Equivalently, the determinant section vanishes on $X^c$ with multiplicity one. On $\PP^1$ a local model for the minimal jump in holomorphic structure  for a  rank 2 bundle $E$ of degree $0$ is given by the one-parameter family of extensions  in $H^1(\PP^1,{\mathcal O}(-2))\cong \C$. If $[e]$ is a generator then the extension $t[e]$ defines a vector bundle $E$:
 $$0\rightarrow {\mathcal O}(-1)\rightarrow E\rightarrow {\mathcal O}(1)\rightarrow 0$$
which is  trivial  if $t\ne 0$ and   ${\mathcal O}(-1)\oplus {\mathcal O}(1)$ if $t=0$.
 
 We can describe such an extension by a transition matrix defined on $\C^*$ by 
 $$g_{10}(t,z)= \begin{pmatrix}\zeta & t \\
0& \zeta^{-1}
\end{pmatrix}$$
and  a global section is given by holomorphic vector-valued functions  $v_0(\zeta), v_1(\tilde\zeta)$ satisfying $v_1(\tilde\zeta)=g_{10}(t,z)v_0(\zeta)$ where $\tilde\zeta=\zeta^{-1}$. The 2-dimensional space of such sections is defined by $v_0(\zeta)=(-ta_1, a_0+a_1\zeta)$. Since $\det g_{10}=1$ it preserves the standard skew form $\epsilon$ on $\C^2$ and taking a basis  $s_1=(t,-\zeta), s_2=(0,1)$ this gives $\epsilon(s_1,s_2)=t$. In a similar fashion, the 4-dimensional space of global sections of $T_F$, with transition matrix $\zeta^{-1}g_{10}(t,z)$ is  given by $v_0(\zeta)=(b-tc_3\zeta, c_0+c_1\zeta+c_2\zeta^2)$.

Applying this to the bundle $T_F(-1)$ together with the skew form defined by $\varpi$.
 we can implement the isomorphism (\ref{tensor}) to  obtain the $\zeta$-dependent skew form on $T_m$
$$\omega_{\zeta}=db\wedge dc_0+(db\wedge dc_1+tdc_0\wedge dc_2)\zeta+(db+tdc_1)\wedge dc_2\zeta^2.$$
When $t=0$ it is well defined and all coefficients are divisible by $db$ -- the characteristic property of an $\alpha$-fold. More importantly, identifying this with 
$(\omega_2+i\omega_3)+2\zeta\omega_1-\zeta^2(\omega_2-i\omega_3)$ we see that $2\omega_2^2=(\omega_2+i\omega_3)\wedge(\omega_2-i\omega_3)=tdb\wedge dc_0\wedge dc_1\wedge dc_2$ and so vanishes transversally on $X^c$.

\noindent 6. Consider now the real structure, given by the antiholomorphic involution $\sigma$, and the real lines parametrised by $M\subset M^c$ . Then $\omega_{\zeta}$ is real and the reality condition on  the coefficients means that the $\omega_i$ are real forms on $M$. As we have seen, $\omega_i^2$ vanishes transversally on $X\subset M$. Moreover, since $\omega_{\zeta}$ is preserved by the holomorphic involution $\tau$ which covers $\zeta\mapsto -\zeta$, we have $\tau^*\omega_1=-\omega_1$ and $\tau^*(\omega_2+i\omega_3)=(\omega_2+i\omega_3)$. We therefore have an $\alpha$-fold $X$ on $M$, which is the fixed point set of an involution. 
\end{prf}

\begin{rmks} 

\noindent 1.  Consider the  real twistor lines  given by $M\backslash X$. The standard twistor approach gives that  the intersection with any fibre $Z_{\xi}$ of $Z\rightarrow \PP^1$  is a local diffeomorphism, and we obtain
a description of the hyperk\"ahler metric as   a $C^{\infty}$ product $\PP^1\times M$.  The twisted holomorphic 2-form is then written as $(\omega_2+i\omega_3)+2\zeta\omega_1-\zeta^2(\omega_2-i\omega_3)$. 

Take $\xi=0$. This fibre is fixed by $\tau$ so $L_m$ and $\tau(L_{m})=L_{\tau(m)}$ meet it at the same point. The restriction map $M\rightarrow Z_0$ therefore factors through the involution $\tau$. Pulling back $\varpi$ gives a closed complex 2-form $\omega_2+i\omega_3$ with $\tau^*(\omega_2+i\omega_3)=(\omega_2+i\omega_3)$ and by construction its restriction to $X$ is  $(\eta_1+i\eta_2)\wedge\varphi$ which is rank 2.  So this is folding in a quite concrete sense: the map from $M$ to $Z_0$ is a folding map, and its image is a manifold with boundary diffeomorphic to $X$.

\noindent 2. It is natural to ask about the geometry on $Z_0$ on the other side of the hypersurface. In fact, composing the real structure on $Z$ with $\tau$ gives a new real structure covering $\zeta\mapsto 1/\bar \zeta$ on $\PP^1$, the reflection in an equator. With this structure $\omega_1$ is imaginary and we have the relation $\omega_2^2=\omega_3^2=-\omega_1^2$ for the three closed forms. The twistor lines which are real for this structure define a 4-manifold with hypersymplectic structure (see e.g.\cite{Hit01},\cite{DS}). In four dimensions this is a Ricci-flat anti-self-dual Einstein metric  of signature $(2,2)$. The analytic continuation of the canonical model to the {\it exterior} of the disc bundle is an example. Inside $M^c$ we thus have two real submanifolds, intersecting in $X$ but mapping to opposite sides of the corresponding hypersurface in the fibre $Z_0$. This is analogous to the real and imaginary axes in $\C$ and the map $z\mapsto z^2$ to $\R$.

\noindent 3. To reverse the construction, if a twistor section has normal bundle $T_F\cong {\mathcal O}(2)\oplus{\mathcal O}$ then the distinguished subbundle ${\mathcal O}(2)$ lifts it canonically to the 4-manifold $\PP(T_F)\rightarrow Z$. It has trivial normal bundle there and so locally this space is a product $ \PP^1\times X^c$.

\noindent 4. The intrinsic differential geometry of a 3-manifold with the boundary data of the theorem has been studied by R.Bryant \cite{Bry1} who  proved elsewhere a hyperk\"ahler extension theorem analogous to the above but for $\beta$-folds \cite{Bry2}.
 
\end{rmks}

\section{Global invariants}

\subsection{Topological invariants}
Our conjectural folded hyperk\"ahler metrics are deformations of the canonical model and so are defined on the 4-manifold $M$ which is the $S^2$-bundle $\PP(1\oplus K)$ over the surface $\Sigma$. The second homology is generated by the classes of a fibre and the zero section in $K\subset \PP(1\oplus K)$. We have three closed 2-forms $\omega_1,\omega_2,\omega_3$ which define de Rham classes in $H^2(M,\R)$. But we also have the involution $\tau$ and in particular $\omega_2,\omega_3$ are invariant by $\tau$.  Since $\tau$ changes the orientation of each fibre it follows that the classes $[\omega_2],[\omega_3]$ evaluated on a  fibre are zero. But also, $\omega_2+i\omega_3$ restricted to $K$ is holomorphic and so vanishes on the zero section. Thus  $[\omega_2]=0=[\omega_3]$. 

From the $SU(\infty)$ point of view $\omega_1$ restricts to the standard symplectic form on a fibre and evaluates to  $4\pi$. From the explicit formula (\ref{om1}) for the canonical model $\omega_1=d(u_3(dx-yd\phi))$ and so on the zero section $x_3=yu_3=1$, $\omega_1=dx\wedge dy/y^2$ and evaluating the integral gives the area of the hyperbolic metric of curvature $-1$ which is $4\pi(g-1)$. 

\subsection{Invariant polynomials}

One of the fundamental features of the moduli space of Higgs bundles is the associated integrable system, based on  evaluating an invariant polynomial of degree $m$ on $\Phi$ to give a holomorphic section of $K^m$ on $\Sigma$. As we remarked above,   we can replace the polynomial $\tr A^m$ on ${\lie su}(n)$ by 
\begin{equation}
p_m(f)=\int_{S^2}f^m\omega
\label{inv}
\end{equation}
for $SU(\infty)$, and this we can do for the Lie algebra and its complexification.

Then for the Higgs field $\Phi$, a section of $p^*K$ on the sphere bundle $M^4$, integration over the fibres of $\Phi^m$ defines likewise a section $\alpha_m$ of $K^m$ on $\Sigma$. This is again holomorphic, for 
 from (\ref{higgs2}) we have $2(\phi_1+i\phi_2)_{\bar z}+\{a_1+ia_2,\phi_1+i\phi_2\}=0$. Then writing $\psi=\phi_1+i\phi_2, a=a_1+ia_2$ we have 
 $$\frac{\partial}{\partial \bar z}\left(\int_{M/\Sigma}\psi^m\omega\right) dz^m=-\frac{1}{2}\left(\int_{M/\Sigma}m\{a,\psi \}\psi^{m-1}\omega\right) dz^m=-\frac{1}{2}\left(\int_{M/\Sigma}{\mathcal L}_{X_a}(\psi^{m}\omega)\right) dz^m.$$
where $X_a$ is the Hamiltonian vector field of $a$. But this integral is zero by Stokes's theorem as ${\mathcal L}_{X_a}(\psi^{m}\omega)= d(i_{X_a}(\psi^{m}\omega))$. This holds for any $SU(\infty)$-Higgs bundle. For the ones we are considering  where the connection lies in $SO(\infty)$ and the Higgs field is symmetric under the involution, the differential $\alpha_m$ is twice the integral over the disc bundle. When $\Phi$ maps the closed disc bundle diffeomorphically to a submanifold $D$ of the cotangent bundle then, using standard local coordinates $(z,w)\mapsto wdz$ this integral is obtained from the canonical holomorphic 1-form $\theta=wdz$:
$$\alpha_m=2\left(\int_{D/\Sigma}w^m\Phi_*\omega\right) dz^m=2\left(\int_{D/\Sigma}w^m\omega_1\right) dz^m$$
using the K\"ahler form $\omega_1$.
\begin{ex} In the example relating to Teichm\"uller space above we can perform this calculation using (\ref{phiham}) to get
$$\alpha_m=2\int_{0\le \theta\le 2\pi}\int_{0\le r \le1}\frac{1}{2^m}(k^{1/2}re^{-i\theta} +ak^{-1/2}re^{i\theta})^m\frac{rdrd\theta}{2\sqrt{1-r^2}} dz^m$$
which is zero if $m$ is odd and if $m=2\ell$ is 
$$\left(\frac{\pi}{2^{2\ell}}{2\ell \choose \ell}\right)^2a^{\ell}dz^{2\ell}$$
and we recover a power of the quadratic differential $q=adz^2$ used in the definition of the Higgs bundle. For the canonical model $q=0$ and all these invariants therefore vanish. 
\end{ex} 

\begin{rmk}
One may ask what  happened to $\alpha_0$ and $\alpha_1$. But $\alpha_0$ is just the area of the $2$-sphere which is the standard $4\pi$. As for $\alpha_1$, in the finite-dimensional case this is set to zero because we are considering the group $SU(n)$ and the trace of the Higgs field is zero. For our case, a translation $wdz\mapsto wdz+a(z)dz$ by a holomorphic section of $K$ preserves the holomorphic symplectic form and takes one solution to another. Setting $\alpha_1=0$ removes this trivial deformation.
\end{rmk}
 \section{First order deformations}
 If a generalization of the higher Teichm\"uller spaces exists then the holomorphic differentials $\alpha_m$ should uniquely determine a folded hyperk\"ahler metric. It makes sense then to look for deformations, and initially infinitesimal deformations, of the canonical model to be determined by a holomorphic differential. The Teichm\"uller example offers a test: we have a genuine deformation in the direction of a quadratic differential.
 
 Differentiating equation \ref{qhiggs} for the quadratic differential $tq$ with respect to $t$  and setting $t=0$ we get $\dot F=2\dot\omega$. Now  the conformal structure is unchanged so $\dot\omega=f\omega$ for some function $f$ But then $\dot F= dd^cf$ and  and so   $dd^cf=f\omega$ and this  implies $f=0$ since $dd^c$ is a negative operator. The infinitesimal variation of the metric is therefore zero and so $\dot\omega_1=0$ for the hyperk\"ahler metric.  From (\ref{phiham})
$$\dot \Phi=\dot w dz= ay^2\bar wdz$$
using the hyperbolic metric $k=y^{-2}$ and $\dot k=0$. Then 
\begin{equation}
\dot \omega_2+i\dot\omega_3=d(\dot w dz)=d(ay^2\bar w dz).
\label{Lomega}
\end{equation}
Introduce the complex vector field  (tangential to  the fibres)
$$X^c=ay^2\bar w\frac{\partial}{\partial w}$$
 and let $X$ be the real part, then we can write (\ref{Lomega}) as  $$\dot \omega_2+i\dot\omega_3=d{\mathcal L}_X wdz={\mathcal L}_X(\omega_2+i\omega_3).$$

 We can generalise this infinitesimal deformation by taking a holomorphic section $\alpha_m$ of $K^m$ and writing it in  local coordinates as $\alpha_m=a(z)dz^m$.  Then define a complex vector field $$X^c=ay^{2m-2}\bar w^{m-1}\frac{\partial}{\partial w}.$$
This is globally well-defined because we have the hermitian form $h=dzd\bar z/y^2$, the canonical 1-form $\theta=wdz$ and the canonical holomorphic Poisson tensor $\pi=\partial/\partial w\wedge \partial/\partial z$. The vector field is then 
$X^c=\pi(\alpha h^{-(m-1)}\bar\theta^{m-1}).$

\begin{prp} If $X$ is the real part of $X^c$, then the  three closed 2-forms ${\mathcal L}_X\omega_i$ are anti-self-dual.
\end{prp}

\begin{prf} The K\"ahler forms $\omega_i$ span the space of self-dual 2-forms at each point and so we need to prove that 
${\mathcal L}_X\omega_i\wedge\omega_j=0$ for all $i,j$.

First consider ${\mathcal L}_X(\omega_2+i\omega_3)={\mathcal L}_Xd(wdz)$. This is 
$d(a\bar w^{m-1}y^{2m-2}dz)$ and since 
 $a(z)$ is holomorphic and $y=(z-\bar z)/2i$ we obtain 
\begin{equation}
(m-1)ay^{2m-3}\bar w^{m-2}(yd\bar w\wedge dz+i\bar wd\bar z\wedge dz).
\label{lie}
\end{equation}
So clearly
$${\mathcal L}_X(\omega_2+i\omega_3)\wedge dw\wedge dz=0={\mathcal L}_X(\omega_2+i\omega_3)\wedge d\bar w\wedge d\bar z$$
and we have 
${\mathcal L}_X\omega_2\wedge\omega_2=0={\mathcal L}_X\omega_3\wedge\omega_2={\mathcal L}_X\omega_2\wedge\omega_3={\mathcal L}_X\omega_3\wedge\omega_3.$
In particular ${\mathcal L}_X(\omega_2^2)=0$ but since $\omega_2^2=\omega_1^2$ we also have ${\mathcal L}_X\omega_1\wedge\omega_1=0$.

It remains to prove that ${\mathcal L}_X(\omega_2+i\omega_3)\wedge\omega_1=0$ for then we shall have ${\mathcal L}_X\omega_2\wedge \omega_i=0$ for $i=1,2,3$ and similarly for $\omega_3$. Taking the Lie derivative of  $\omega_3\wedge \omega_1=0=\omega_2\wedge \omega_1$ we shall then have the same result for ${\mathcal L}_X\omega_1$.

Recall then that $$\omega_1=d(u_3(dx-yd\phi))=du_3\wedge (dx-yd\phi)-u_3dy\wedge d\phi$$ and 
$u_3y=\sqrt{1-y^2\vert w\vert^2}, e^{2i\phi}=w/\bar w.$
Since from (\ref{lie}) $dz=d(x+iy)$ is a factor of ${\mathcal L}_X(dw\wedge dz)$ we calculate first 
\begin{equation}
\omega_1\wedge(dx+idy)=idu_3\wedge dx\wedge dy-u_3d\phi\wedge dx\wedge dy-ydu_3\wedge d\phi\wedge dx-iydu_3\wedge d\phi\wedge dy
\label{wedge}
\end{equation}
and now, from the form of ${\mathcal L}_X(dw\wedge dz)$ in (\ref{lie}), take the exterior product with $yd\bar w+i\bar w d\bar z$. 

 Using 
$$2id\phi=\frac{dw}{w}-\frac{d\bar w}{\bar w},\quad du_3=-(u_3^2+\vert w\vert^2)\frac{dy}{yu_3}-\frac{1}{2u_3}(\bar wdw+wd\bar w)$$ in the expression (\ref{wedge}) we get zero. 

\end{prf}

Setting $\dot\omega_1=0, \dot\omega_2= {\mathcal L}_X\omega_2$ and $\dot\omega_3= {\mathcal L}_X\omega_3$ from the proposition we certainly have a first order deformation of the algebraic equations
$$\omega_1^2=\omega_2^2=\omega_3^2=0 \qquad \omega_1\wedge \omega_2=\omega_2\wedge \omega_3=\omega_3\wedge \omega_1=0$$
by closed forms. We also have a deformation which is orthogonal to just rotating the forms $\omega_i$. 

Infinitesimal  deformations of Ricci-flat metrics define elements in the null space of the Lichnerowicz Laplacian acting on trace-free symmetric tensors. On a hyperk\"ahler 4-manifold these are constructed from tensor products of the covariant constant self-dual 2-forms $\omega_i$ and closed anti-self dual 2-forms, thus the proposition gives us potential deformations of hyperk\"ahler metrics (and in the compact case would give actual deformations). 

We can calculate the infinitesimal variation in the global invariants using the above. Since $\dot\omega_1=0$ the cohomology class is unchanged. 
Since  $\omega_1^2$ is fixed to first order then so is the fold, so we only have to vary  $w$ which is the vector field $X$ applied to $w$. The variation in the moment of $w^k$ is therefore given by

 $$\int_{y^2\vert w\vert^2<1}kw^{k-1}a\bar w^{m-1}y^{2m-2}y^2\frac{dwd\bar w}{\sqrt{1-y^2\vert w\vert^2}}.$$
The measure is $S^1$-invariant and so the multiple integral vanishes unless $m=k$. In this case put $u=yw$ and we get 
$$a\int_{\vert u\vert<1}\vert u\vert^{2k-2}\frac{dud\bar u}{\sqrt{1-\vert u\vert^2}}=2\pi a\int_0^1r^{2k-2}\frac{rdr}{\sqrt{1-r^2}}=\pi a \frac{4^k (k!)^2}{(2k)!}.$$
This is consistent with our conjecture that the differentials determine the metrics.
\section{$S^1$-invariance}
The usual Teichm\"uller space is embedded in its higher analogue by the homomorphism $SL(2,\C)\rightarrow SL(n,\C)$ given by the irreducible $n$-dimensional representation. Restricted to the compact real form we take the $SU(n)$-Higgs bundle  induced from the $SU(2)$-Higgs bundled discussed in Section 6.2 and restricted to the split form we take the positive hyperbolic $SL(n,\R)$-representation coming from the uniformizing representations. The scalar $S^1$-action $\Phi\mapsto e^{i\phi}\Phi$ preserves the higher Teichm\"uller spaces but acts non-trivially on all differentials except zero and so the only solution fixed by this action is the canonical solution. If our conjecture holds, then the canonical model should be the only $S^1$-invariant folded hyperk\"ahler metric on a  disc bundle in the cotangent bundle. Here, where we regard the Higgs field as the canonical 1-form, the action is just scalar multiplication by a unit complex number in the fibre of $T^*\Sigma$. 

This circle action acts on $\omega_2+i\omega_3=d(wdz)$ by $e^{i\phi}$ and preserves $\omega_1$ and for this type of action in four dimensions there is a well-known differential equation called either the Boyer-Finley equation \cite{BF} or the $SU(\infty)$-Toda equation \cite{Tod}.

Locally the Ansatz for the metric is 
\begin{equation}
u_t(e^u(dx^2+dy^2)+dt^2)+u_t^{-1}(d\tau+u_ydx-u_xdy)^2
\label{Boy}
\end{equation}
and for the K\"ahler form 
\begin{equation}
\omega_1= u_te^udx\wedge dy+dt\wedge (d\tau+u_ydx-u_xdy).
\label{kform}
\end{equation}
Requiring this to be closed gives the Boyer-Finley  equation  
\begin{equation}
u_{xx}+u_{yy}+(e^u)_{tt}=0.
\label{toda}
\end{equation}
The holomorphic symplectic form is 
$\omega_2+i\omega_3=d(wdz)$
where $w=\sqrt{2} e^{u/2+i\tau}$. 
\begin{ex}
For the canonical model $u=\log(1-t^2)-2\log y$, giving $u_x=0$ and $u_y=-2/y$. Hence $u_{yy}=2/y^2$ and $(e^u)_{tt}=-2/y^2$. Then $t=0$ is the fold and at $t=-1$, $u_te^u=2/y^2$ which gives the hyperbolic metric $2y^{-2}(dx^2+dy^2)$ on the zero section of the cotangent bundle. 
\end{ex}

\begin{rmk}
One form of the 
2-dimensional  Toda lattice is the differential-difference equation $$v_{xx}+v_{yy}+e^{v_{n+1}}-2e^{v_n}+e^{v_{n-1}}=0.$$
The limit where the difference operator becomes the second derivative is the reason for the $SU(\infty)$ terminology and is consistent with our use in this paper. 
\end{rmk}

We need to give  the local formulae above a more geometric interpretation. The Killing field is $\partial/\partial\tau$ so from (\ref{kform}) the moment map with respect to $\omega_1$ is  $-t$. Setting $t$ to be a constant and taking the quotient by the circle action gives the K\"ahler quotient as the metric  $4u_te^u(dx^2+dy^2)$ so  $z=x+iy$ is a local complex coordinate on the quotient. The function $u$ itself is defined by $w\bar w=2e^u$. 

If we focus now on the $t$-dependent metric $g=w\bar w dzd\bar z=2e^u(dx^2+dy^2)$ on a surface then its Gaussian curvature $\kappa=-e^{-u}(u_{xx}+u_{yy})/4$ and then the Toda equation (\ref{toda}) may be written in the more geometrical form 
\begin{equation}
g_{tt}=4\kappa g
\label{toda1}
\end{equation} 
for a $t$-dependent family of metrics on a compact surface.

In our setting $wdz$ is the canonical 1-form on $T^*\Sigma$ restricted to a disc bundle. The $S^1$-invariance $w\mapsto e^{i\phi}w$ means our boundary value data for the fold is the unit cotangent bundle for a  hermitian metric .  

In the Higgs bundle picture $\omega_1$ restricted to the fibre is the restriction of the standard symplectic form on $S^2$ to the disc  $x_3\ge 0$ and $x_3$ itself is the moment map for the circle action so $t=-x_3$ and with this normalisation the fold is $t=0$ and the zero section $t=-1$ as in the canonical model.  Note from formula (\ref{kform}) that integrating $\omega_1$ over each disc is integrating $dt\wedge d\tau$ for $-1\le t\le 0,0\le \tau\le 2\pi$ which gives $2\pi$, the correct normalisation.

For regularity we assume smoothness of the forms $\omega_i$ on the $S^2$-bundle $\PP(1\oplus K)$ and on this space the Higgs field is described locally as $\phi_1dx+\phi_2dy$ where $\phi_1,\phi_2$ are even functions. The metric $g=(\phi_1^2+\phi_2^2)(dx^2+dy^2)$ and since the $\phi_i$ are even $\partial \phi_i/\partial x_3=0$ on the fold. 
Hence $g_t=0$ at the fold. We now have boundary conditions $g_t(0)=0$ and $g(-1)=0$ for the equation (\ref{toda1}). In fact, we are interested in a solution over the whole sphere bundle, invariant under an involution, so we could  take equivalently the boundary conditions $g(-1)=0=g(1)$. 
 
 \begin{thm} The only $S^1$-invariant hyperk\"ahler metric on $\PP(1\oplus K)\rightarrow \Sigma$ which has an $\alpha$-fold on a circle bundle is the canonical model for a hyperbolic metric on $\Sigma$.
 \end{thm}
 \begin{prf}
 
\noindent 1.  Note that since $g=w\bar wdzd\bar z$ the conformal structure is independent of $t$, so if $g=fh$ for some fixed metric $h$, the area form $v_g=fv_h$. Now integrate (\ref{toda1}) against $v_h$ and use Gauss-Bonnet.  If $A(t)$ is the area of the metric $g(t)$ on $\Sigma$ we get 
$$A_{tt}=16\pi(1-p)$$
where $p$ is the genus, and hence $A(t)=(a+bt+ct^2)$ for constants $a,b,c$. But the metric vanishes at $t=\pm1$ so $b=0=a+c$ and we have 
$$A(t)=8\pi (t^2-1)(1-p).$$
and  since $-1\le t\le 1$ the genus must be greater than 1. 

\noindent 2. Consider  the metric $h=g/2(1-t^2)$, then  the area is fixed and the metric is regular at $t=\pm 1$. 
 The right hand side of the equation (\ref{toda1}) is the Gauss-Bonnet integrand which is unchanged under constant rescaling so the equation for $h$ is
$$(1-t^2)h_{tt}-4th_t-2h=4\kappa h$$
and clearly $h$ of constant curvature $-1/2$ satisfies this.  These metrics are all conformally related so we can write $h=fg_{H}$ where $g_{H}$ is the hyperbolic metric and then we get 
$(1-t^2)f_{tt}-4tf_t=2\Delta_H\log f $
 where $\Delta_H$ is the Laplace-Beltrami operator for the hyperbolic metric. 
 
 Rewrite this as 
 $((1-t^2)^2f_t)_t=2(1-t^2)(\Delta_H\log f). $
 Integrating by parts and using the fact that $f$ is finite at $t=\pm1$, we have 
$$2\int_{-1}^1(1-t^2)f\Delta_H\log f dt =\int_{-1}^1((1-t^2)^2f_t)_tfdt=-\int_{-1}^1(1-t^2)^2f^2_tdt.$$
But
$f\Delta_H\log f=\Delta_Hf+f^{-1}(df,df)$ so integrating over $[-1,1]\times \Sigma$ we get 
$$-\int_{-1}^1(1-t^2)^2\int_{\Sigma}f^2_tdt\,v_H=2\int_{-1}^1(1-t^2)\int_{\Sigma}(\Delta_Hf+f^{-1}(df,df))dt\,v_H$$
which, since the integral of $\Delta_Hf$ over $\Sigma$ vanishes, is a contradiction unless $f$ is a constant on $[-1,1]\times \Sigma$.  
\end{prf}

\begin{rmks} 1. There do exist $S^1$-invariant solutions to the Higgs bundle equations other than the canonical one. The holomorphic Higgs bundle is of the form
$$E=L\oplus L^* \qquad \Phi= \begin{pmatrix}0 & s \\
0& 0
\end{pmatrix}$$
where $s$ is a holomorphic section of $L^2K$. But when $L\ne K^{-1/2}$, $s$ has zeros and then the Higgs field itself vanishes, which gives a rather different singularity for the hyperk\"ahler metric, apart from the circle bundle fold.  

2. The uniqueness result above, coupled with the local existence theorem of Section 7,  shows that the unit cotangent bundle of  a Riemannian metric which has a  global folded hyperk\"ahler extension is a hyperbolic metric. If our conjecture holds, then the general situation must involve a geometry of Finsler type. 
\end{rmks}

 \end{document}